# On the General Shape of Scales on Slide Rules

*István Szalkai*[*]

ver: 2017.06.10.,

**Abstract**

In this article we continue our investigations started in our previous papers [2] and [3]: we give a thoroughful explanation of the general properties ("shape" and others) of different, general scales, corresponding to different (all possible) mathematical functions $f(x)$. Though we do not plan here any new scales for specific tasks, we mention and analyse many examples. These observations and statements might help us to plan and realize any specific scale we need in the future, or to understand better any interesting or miserable properties of old scales.

Both for illustrations, investigations and planning we suggest our Javascript animation " Rubber Band.html" [6]. However, we suggest the Readers to think on the function and the scale with paper and pencil, as it is shown on Figure 1, and after to check her/himself by [6]. In the last chapter we list some possible applications of our several new scales, different from [2]. Some new scales are visualized in professional way on Ace Hofman's webpage [1].

### 0  How they are made ?

Translating (moving) the laths of the slide rules realize physical addition of lengths. However, as in many places ([2,3,8]) is explained, we usually write other numbers $x,y,z$ at the distances $d_x = f(x)$, $d_y = g(y)$ and $d_z = h(z)$, which helps us to easily read the result

$$h(z) = f(x)+g(y) , \quad \text{i.e.} \quad z = h^{-1}(f(x)+g(y)) , \qquad (1)$$

instead of the addition of the physical distances $d_z = d_x + d_y$, where $h^{-1}$ is the inverse function[**] of the function $h$. For example, the base scales "$x$" (named C and D scales [7],[8]) on usual slide rules we have $d_x = f(x) = log(x)$, $d_y = g(y) = log(y)$ and $d_z = h(z) = log(z)$, which gives us the relation

$$log(z) = log(x)+log(y) = log(x \cdot y) , \quad \text{i.e.} \quad z = x \cdot y . \qquad (2)$$

For practical reasons $f(x)$, $g(y)$ and $h(z)$ must be <u>strictly monotone functions</u> (either increasing or decreasing). When adjusting the size of the scale <u>and</u> the interval of numbers to be printed on the latch we zoom the scale, that is we choose an appropriate (but arbitrary) physical unit $u$ and we really use $d_x = u \cdot f(x)$. In this paper we do not deal with $u$.

A comfortable and illustrative construction of the scale $d_x = f(x)$ is shown on Figure 1. We start with an equidistant sequence of numbers $x$ as $1,2,3,...,10$, we draw vertical lines from $x$ to the graph of the function $f(x)$ and then, from the intersection point we start a horizontal line to the $y$ axis. This horizontal line on the $y$ axis meets us the distance $d_x = f(x)$ for the number $x$, and the Reader can recognize the well known "$log$" scale on the $y$ axis, which is printed on slide rules. (We also constructed the lower scale for the values $x = 0.1, 0.2, ... , 1.0$ on Figure 1.)

For decreasing functions the scale on the $y$ axis is downward directed: the horizontal lines move from top to bottom while the vertical ones from left to right.  $x$ also could be negative, in which case we again draw the vertical lines from the $x$ axis to the graph and the horizontal ones from the graph to the $y$ axis. Keep in mind that only for (strictly) monotone functions is possible to construct a scale. (See also (iii) in Section 2.)

The purpose of the present paper is a thoroughful discussion of the properties of the scales concerning to (almost) all possible (mathematical) functions $f(x)$. Though the above graphical method can be eliminated by computer graphics and printers, it helps us to understand the shape and behaviour the the scale and to plan the length and content of it. We suggest the Readers <u>first</u> to sketch the graph of the function in question and the derived scale on the $y$ axis, and <u>then</u> to read on this article and make computer experiments by [6]. Section 8 contains some example scales with explanations.

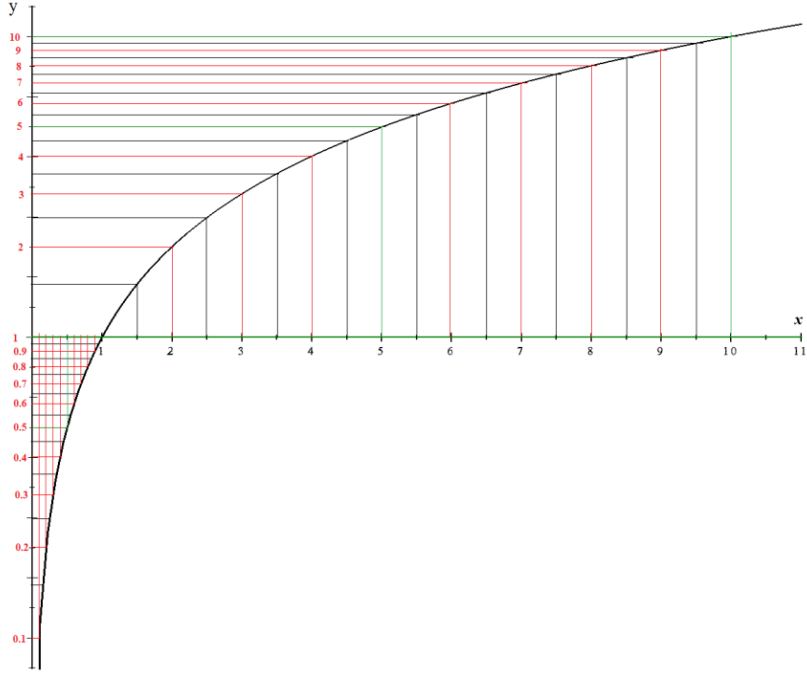

**Figure 1.:** Constructing the *log* scale

## 1 Traditional Scales and Other Examples

    In this article we call the scales after the functions *f(x)* which really gives us the physical distance $d_x=f(x)$. However, the usual naming of the scales on "traditional" slide rules is different, which we are forced to discuss just at the beginning, here. The base scales C and D, named "*x*" use the distance function $d_x = lg(x)$ where *lg* is the logarithm of base 10, these scales are, in fact *lg(x)* scales. The other scales are <u>one variable</u> functions *t(x)* of *x* of scale D, and so they are named "*t(x)*" scales (see [7],[8]). For example, the scale L is marked by " $x^3$ " which means, that after choosing any value *x* on scale D, the value exactly on the same place/distance on scale K (use the hairline) is $x^3$, that is $t(x)=x^3$ on scale K. The title of scale L is "*lg(x)*", so the value *x* on scale D corresponds to $t(x)=lg(x)$ on scale L.

    A closer calculation reveals the real formula of the K,L and other scales. If the distances of the value *x* on scales D and of $y=t(x)$ on X are equal, then we must have $d_x = lg(x) = d_y = \tau(y) = \tau(t(x))$ where $\tau(y)$ is the distance-function for the scale X we are looking for. We easily get $x=t^{-1}(y)$ and so $\tau(y)= lg(t^{-1}(y))$

$$d_y = \tau(y) = lg(t^{-1}(y)) \,. \qquad (3)$$

For example, scale K ($y=x^3$) has the distance-function $\tau_K(y) = lg(\sqrt[3]{y}) = lg(y)/3$, and look: scale K is the third of D. Or, for scale L ($y=lg(x)$) we have $\tau_L(y) = lg(10^y) = y$, and really: L is in fact an equidistant scale, as an ordinary ruler or measuring tape! We calculate some other traditional scales' formulae at the end of this Section.

    Before reading on the next Section, we advice to the Reader to sketch the graph and the scale with paper and pencil for different, usual and unusual functions *f(x)*. For example the following ones: *cx*, $x^2$ and $\sqrt{x}$ for $0\leq x$, $x^3$ and $\sqrt[3]{x}$ for $x\in R$, $\sqrt{1-x^2}$ for $0\leq x\leq 1$, $\sqrt{x^2-1}$ for $1\leq x$, $\sqrt[3]{1-x^3}$ and $\sqrt[3]{x^3-1}$ for $x\in R$, $1/x$ for $0<x$ and for $x\leq 0$ (two drawings), all the functions $x^\alpha$ and $\sqrt[\alpha]{1-x^\alpha}$ and $\sqrt[\alpha]{x^\alpha-1}$ for all possible $\alpha\in R$ and $x\in R$, $sin(x)$ and $sin^2(x)$ for $-\pi/2\leq x\leq +\pi/2$, $cos(x)$ and $cos^2(x)$ for $0\leq x\leq \pi$, $tan(x)$ for $-\pi/2<x<+\pi/2$, $atan(x)=tan^{-1}(x)$ (inverse of *tan(x)*) for $x\in R$, $e^x$, $a^x$ and $1-e^{-\alpha x}$ for $x\in R$ and various $a,\alpha\in R^+$, $log_a(x)$ for $0<x$ and various $a\in R$, $log_x(10) = 1/log_{10}(x)$ either for $0<x<1$ or for $1<x$, $\frac{1}{1+x^2}$ for $0\leq x$, $e^{-mx^2}$ for $0\leq x,m$, $f_G(x) = r \cdot arccos\left(\frac{r}{r+x}\right)$ for $0\leq x$ (see [2],[3]), $log(sin(x))$ for $0<x\leq +\pi/2$, $log(tan(x))$ for $0<x<+\pi/2$, $\frac{-1}{2}log\sqrt{1-x^2/c}$ for $0\leq x<c$ (Lorentz transformation [2]), $sh(x)=sinhyp(x) = \frac{e^x-e^{-x}}{2}$ and $Arsinhyp(x) =sh^{-1}(x)= ln(x + \sqrt{x^2+1})$ for $x\in R$, $ch(x)=coshyp(x)= \frac{e^x+e^{-x}}{2}$ for $x\in R$ and $Arcoshyp(x)=ch^{-1}(x) = ln(x + \sqrt{x^2-1})$ for $0\leq x$, $tanhyp(x)=\frac{e^x-e^{-x}}{e^x+e^{-x}}$ for $x\in R$ and $Artanhyp(x)=tanhyp^{-1}(x)= \frac{1}{2}ln\left(\frac{1-x}{1+x}\right)$ for $-1<x<+1$ (**), $\Phi(x)$ for $x\in R$ ([cumulative distribution function](#) [9],[10], which is similar but not equal to the [error function](#)

[10] of the normal distribution), $\Phi(\ln(x))$ for $x \in R^+$ (lognormal distribution), $f_{921}(x)=x^5-3x^3$ for $-3/\sqrt{5}<x<3/\sqrt{5} \approx$ 2.236, $f_{923}(x)=1-x^2/2+x^4/24$ for $0<x<\sqrt{6}$ ($x_0=\sqrt{6-2\sqrt{3}} \approx 1.593$), $f_{942}(x)=\frac{x^3}{x^2-x-2}$ either for $x<1-\sqrt{7}\approx-1.646$ or for $-1<x<2$ or for $3.646 \approx 1+\sqrt{7}<x$, $f_{953}(x)=x+\sqrt{1-x^2}$ for $-1 \leq x \leq \sqrt{2}/2 \approx 0.7071$ ($x_0=-\sqrt{2}/2 \approx -0.7071$), ..., and so on. Some scales are shown in the last but one Section.

Now we list the real formulae for some scales of traditional slide rules.

**A:** "$x^2$" has $d_x=t(x)=lg(x)/2$,      **C,D:** "$x$" has $d_x=t(x)=lg(x)$,
**CI:** "$1/x$" has $d_x=t(x) = -lg(x)$,      **K:** "$x^3$" has $d_x=t(x)=lg(x)/3$,
**L:** "$log(x)$" has $d_x=t(x)=x$,      **Ln:** "$ln(x)$" has $d_x=t(x)=x*lg(e)$,
**LL3:** "$e^x$" has $d_x=t(x)=lg(ln(x))$,      **LL10:** "$10^x$" has $d_x=t(x)=lg(lg(x))$,
**LL03:** "$e^{-x}$" has $d_x=t(x)=lg(-ln(x))$,      **R1:** "$\sqrt{x}$" has $d_x=t(x)=2*lg(x)$,
**H1:** "$\sqrt{1+(0.1x)^2}$" has $d_x=t(x)=lg(10\sqrt{x^2-1})$,      **H2:** "$\sqrt{1+x^2}$" has $d_x=t(x)=lg(\sqrt{x^2-1})$,
**P1:** "$\sqrt{1-(0.1x)^2}$" has $d_x=t(x)=lg(10\sqrt{1-x^2})$,      **P2:** "$\sqrt{1-x^2}$" has $d_x=t(x)=lg(\sqrt{1-x^2})$,
**S0:** "$sin(x)$" has $d_x=t(x)=lg(arcsin(x))=lg(sin^{-1}(x))$,      **S:** "$asin(x)$" has $d_x=t(x)=lg(sin(x))$,
**T0:** "$tan(x)$" has $d_x=t(x)=lg(arctan(x))=lg(tan^{-1}(x))$,      **T:** "$atan(x)$" has $d_x=t(x)=lg(tan(x))$,
**Sh:** "$Arsinhyp(x)$" has $d_x=t(x)=lg(sinhyp(x))$,      **Ch:** "$Arcoshyp(x)$" has $d_x=t(x)=lg(coshyp(x))$,
**Th:** "$tanhyp(x)$" has $d_x=t(x)=lg(tanhyp(x))$,

**Sh0:** "$sinhyp(x)$" has $d_x=t(x) = lg(Arsinhyp(x)) = lg(sinhyp^{-1}(x)) = lg\left(\ln\left(x+\sqrt{x^2+1}\right)\right)$,

**Ch0:** "$coshyp(x)$" has $d_x=t(x) = lg(Arcoshyp(x)) = lg(coshyp^{-1}(x)) = lg\left(\ln\left(x+\sqrt{x^2-1}\right)\right)$,

**Th0:** "$tanhyp(x)$" has $d_x=t(x) = lg(Artanhyp(x)) = lg(tanhyp^{-1}(x)) = lg\left(\frac{1}{2}\ln\left(\frac{1-x}{1+x}\right)\right)$.

## 2 Basic observations

Below you can find many trivial observations, too, but here we would like to give a complete list of connections between the properties of functions and of the corresponding scales. Well, a thoroughful, complete and general discussion is impossible since there are many different various properties of (infinitely) many functions, moreover most functions have more different mix of properties, too. But we try our best.

The mathematical background for translation and zooming is discussed in [3], [4].

**(o)** We denote the starting point of the scale (where the distance $d=0$) by $S_1$, this is usually but not exclusively the left end of the scale. Further, we denote by $x_0$ the number $x_0$ which mut be printed to the place of $S_1$. i.e. $d_{x_0}=f(x_0)=0$. For example, on the C and D scales $x_0=1$ since $d=log(1)=0$.

Recall, that the value "$x$" is written at (signed geometric) distance $d_x=f(x)$ from the (geometric) point $S_1$. So, the value "$x$" is written to the right / to the left of $S_1$ if and only if $f(x)$ is positive/negative. This implies, that $S_1$ is at the left/right end of the scale if and only if $f(x)$ has only nonnegative/nonpositive values. If $f(x)$ has both poisitive and negative values, then $S_1$ must be inside of the scale!

For example, the function $f(x)=lg(sin(x))$ is monotone increasing and all values of it are negative, so on this scale $S_1$ is at the right end and $x_0=\pi/2$ since $f(\pi/2)=0$. The scale CI ("$1/x$"), where $f(x)=-lg(x)$ has starting point $S_1$ at his right end. Other examples are $f_1(x)=\sqrt[3]{x^3-1}$ and $f_2(x)=\sqrt[3]{1-x^3}$ where $S_1$ is inside of the scale (both $f_1(x)$ and $f_2(x)$ has both negative and positive values), and $x_0=1$.

Unfortunately the most important question "*where is the beginning* ($S_1$) *of the scale and what value $x_0$ has to be written to it?*" in general is very complicated, so it needs a separate Section, numbered by 3. Have in mind, that $S_1$ and $x_0$ always must be present on the scale!

What about zooming? No problem in the 21th century: the computer screen (e.g. in [1] or [6]) is better than a magnifying glass on the window of mechanical slide rules! Try out $10^x$, $arctan(x)+\pi$ or $\Phi(x)$ : $S_1$ is at the mark "$-\infty$" which must be always on the scale, which makes zooming difficult. We discuss this problem in more detail at the end of Section 3.

When translating the function $f(x)$ vertically to $g(x)=f(x)+v$ ($v \in R$ is any real number), what the difference between the scales of $f(x)$ and of $g(x)$ will be? Nothing for the first glance: just translate (move) the slide with the scale on it by the distance $v$ (more precisely: $v$ units). However, moving the slide, $S_1$ moves together all the marks, but for $g(x)$ we must have a new $S'_1$. We discuss this problem in more detail at the end of Section 3, untill then compare the scales for $1/x$ to $1+1/x$ or $arctan(x)$ to $arctan(x)+\pi/2$ .

Another trivial but important point is, that (in most cases) the scale contains only an interval (a subset) of the whole Dom(*f*), the domain of *f(x)*. This means, that all the below criteria should be required only to that interval. For example, many functions (like $x \cdot ln(x)$) have both in- and decreasing parts, while only monotonic functions can be drawn on a scale (see (iv)). No problem: we draw one monotone increasing interval on one scale, and the other, decreasing interval on the other scale.

**(i)** The sign of *f(x)* is merely different to the sign of *x* itself! Why not write negative numbers "*x*" (with their real sign) on the scales? For addition real numbers in elementary schools or in geometry the scale *f(x)=x* is the best one, since it illustrates both negative and positive numbers "*x*". You can also become familiar with scales $f(x)=x^3$, $f(x)=\sqrt[3]{x^3-1}$ or $f_{953}(x)=x+\sqrt{1-x^2}$ (-1≤x≤√2/2). Even, for the function $f(x)=x^2$ we could write (or at least imagine) both the "+*x*" and "-*x*" values to the same place, but this might be confusing a bit (see also (iii)).

**(ii)** The (part of the) function set on a scale must be continuous (more precisely: must be twice continuously derivable). This is for practical reason: we can not write the corresponding value "*x*" to each geometric point, so in most cases we have to approximate (interpolate) the missing values. This is a bad news for drawing discrete cumulative distribution functions in probability theory. Though the slide rule for binomial coefficients $\binom{n}{k}$ [11] seemly contradicts to this remark, but that slide rule is not based on a single function *f(x)*.

**(iii)** The function drawn on a scale must be either strictly monotone increasing or strictly monotone decreasing. (In the language of calculus: its derivative must have constant sign: either nonnegative or nonpositive.) As we mentioned in (i): if a function has both type intervals, then these intervals must be drawn to different scales. See e.g. $f(x)=x \cdot log(x)$ for 0<*x*<1/e and for 1/e<*x* .

If a function is constant in an interval [*a,b*] (i.e. *f(x)* is not strictly monotonic on [*a,b*]), then all the numbers between *a* and *b* must be positioned to the same geometric pont, which must be clearly indicated on the scale.

**(iv)** The numbers written on the scale are increasing from left to right or from right to left exactly in the case the function *f(x)* is monotone increasing or decreasing. (In calculus language: its derivative is nonnegative or nonpositive.) See e.g. $x^3$ and $1/x$ .

**(v)** The scale is exactly equidistant just in the case *f(x)=cx+b* is a linear function for some constants *c,b*∈R, *c*≠0. See for example scale L ("*log(x)*") on traditional slide rules. For slanting asymptotes see also (xiii).

**(vi)** The numbers (or ticks instead them) are rare or dense if and only if the function *f(x)* is gently or steeply sloping (its derivative is small, i.e. close to 0, or has large absolute value). Compare for example the functions $f(x)=x^3$ at points 0 and 2, $f(x)=\sqrt[3]{x^3-1}$ at points 0 and 1, or *f(x)=log(x)* at points 0.5 and 2 on Figure 1. The scale for the function $f_{923}(x)=1-x^2/2+x^4/24$ for $0<x<\sqrt{6}$ is also worth seeing ($x_0=\sqrt{6-2\sqrt{3}}$≈1.593): the rate of its derivative keeps changing from point to point. See also (viii).

**(vii)** The derivative *f* ' might be zero (the function has a horizontal tangent line), which means that *f(x)* is almost constant, but this phenomenon is allowed only some separated points, like $f(x)=x^3$ or $f(x)=\sqrt[3]{1-x^3}$ at the point *x*=0. By (ii), the values of the numbers *x* (tick for the same magnitude) are very dense around these points. Of course, the range of the density depends on the magnitude of the derivative *f* ' : compare e.g. the functions $f_1(x)=x^2$ and $f_2(x)=x^{10}$ around the points *x*=1 and *x*=2.

Clearly *f* ' can not be 0 on any subinterval [$x_1,x_2$], since in this case *f* would be constant between $x_1$ and $x_2$, this case was discussed in (iii).

**(viii)** As *f* ' and *f* " (the first and second derivatives) are different, the rate of the denseness and the changing of it are different, too. The functions e.g. $x^3$, $\sqrt[3]{1-x^3}$ or *arctan(x)* illustrate well this relation.

**(ix)** The (physical) length of the scale depends on Im(*f*) (image or range of *f* , which is the set of values of *f*), since we write the numbers *x* (for the choosen domain *a*≤*x*≤*b*) at the (physical) distances *f(x)*. Well, when we plan or zoom our scale (e.g.using [6]), we modify *a* and *b*, but the (numerical) content which must fit in on the scale is the interval [*f(a), f(b)*], in fact. For example, when constructing the *tan(x)* scale, we can not choose the whole interval (-π/2, +π/2) since it would require an infinite length scale, so we can choose only a sub-interval [*a,b*]⊂(-π/2, +π/2) ! The properties mentioned in (vi) through (viii) are also important when we choose the interval [*a,b*].

As we mentioned in (o), $S_1$ must always be present on the scale, so $x_0$ must always be included in the set [*a,b*]. (For determining $x_0$ see the next Section.)

The functions *sin(x)* and $\sqrt{1-x^2}$ are useful examples for observing and understanding these phenomena. No problem to draw the scales for many other functions, for example *log(x)*, $x^2$ or $x^3$ : we can choose any domain [*a,b*] ⊂R$^+$ as long as we have a latch (or paper slip) long enough for containing the values in [*f(a), f(b)*].

Other functions, like $atan(x)=tan^{-1}(x)$, $tanhyp(x)=\frac{e^x-e^{-x}}{e^x+e^{-x}}$ or $\Phi(x)$ provide us another problem. Though their range is a finite interval $((-\pi/2,+\pi/2)$, $(-1,+1)$, $(0,1)$ resp.), but the whole set of real numbers (of infinite length) must be fitted on our scale of finite length, since these functions' domain is the whole wet R. In this case we again have to restrict ourselves to an appropriate subinterval $[a,b]\subset R$.

**(x)** In the case $lim_{x\to\infty} f(x) = d$ all the real numbers $x$, larger than a specific $x_1$, are crowded on the scale at the distance $d$ from $S_1$ (since this limit means "$|f(x)-d|<\varepsilon$ for any $\varepsilon>0$ and $x>x_1$"). See for example $f(x)=1/x$ $(d=0)$, $f(x)=atan(x)$ $(d=\pi/2)$, $f(x)=tanhyp(x)$ $(d=1)$, $f_G(x)= r \cdot arccos\left(\frac{r}{r+x}\right)$ $(d=r)$ or $\Phi(x)$ $(d=1)$.

**(xi)** In the case $lim_{x\to c} f(x) = \infty$ we can not represent the value $x=c$ on the scale, since its distance from $S_1$ would be <u>infinitely large</u>. Moreover, the values close to $c$ are also far enough not to be able to be represented on our scale! See for example $f(x)=1/x$ $(c=0)$, $f(x)=tan(x)$ $(c=\pi/2)$, $f(x)=Artanhyp(x)=tanh^{-1}(x)=\frac{1}{2}log\left(\frac{1-x}{1+x}\right)$ $(c=-1$ and $c=+1)$.

**(xii)** In the case $lim_{x\to\infty} f(x) = \infty$ we can choose any interval $[x_0,b]\subset R^+$ as long as we have a scale long enough for containing the values in $[f(x_0), f(b)]$, as we discussed in (ix).

**(xiii)** If $f(x)$ has a <u>slanting asymptote</u>[***], i.e. $f(x)\approx cx+b$ for some $c,b\in R$, $c\neq 0$, then the scale is nearly equidistant for $x>x_1$ for some $x_1\in R$, since $f(x)$ is almost equal the straight line (linear function) $y=cx+b$, and equidistant scales were discussed in (v). Have a look at the functions e.g. $\sqrt[3]{x^3-1}$ (see Figure 3), $\sqrt[3]{1-x^3}$, $\sqrt{x^2-1}$ or $x+1/x$. The Reader might be curious also of the function $f_{942}(x)=\frac{x^3}{x^2-x-2}=x+1+\frac{3x+2}{x^2-x-2}$ and its scale either for $x<1-\sqrt{7}\approx -1.646$ or $3.646\approx 1+\sqrt{7}<x$. $f_{942}(x)$ is strictly monotone increasing on these intervals and has the asymptote $y=x+1$ both for $x\to -\infty$ and $x\to +\infty$. $S_1$ and $x_0$ for $f_{942}(x)$ are discussed at the end of the next Section.

Horizontal and vertical asymptotes were discussed in (x) and (xi).

In Section 5 we discuss some more special property.

## 3  Where $S_1$ should be and what value $x_0$ is here

$S_1$ is the geometrical (physical) starting point of the scales and $x_0$ the number (or symbol) we write at this position, so these base data are crucial when constructing our new scale. $S_1$ must be denoted on the scale in all circumstances, $x_0$ might be omitted. The other numbers $x\in Dom(f)$ are printed at distances $d_x=f(x)$ from $S_1$ : to the left if $f(x)$ is negative, and to the right if $f(x)$ is positive. Figure 1 clearly explains the role of $S_1$ : our scale is constructed on the $y$ axis by the horizontal lines and $S_1$ is at the 0 point of the $y$ axis, at the origo, which concerns to the horizontal line the $x$ axis. $x_0=1$ and look: the numbers $x<1$ are below $S_1$ since $f(x)<0$ for $x<1$, and the numbers $x<1$ are above $S_1$ since $f(x)>0$ for $x>1$.

If the function has a root $x_0$ as $f(x_0)=0$ then clearly $x_0$ is found, moreover $x_0$ is unique since every strictly monotonic function can have at most one root! In the case $x_0$ does exists, we have three cases: $f(x)>0$ for all $x\in Dom(f)$ implies that $S_1$ must be at the <u>left</u> end of the scale, $f(x)<0$ for all $x\in Dom(f)$ implies that $S_1$ must be at the <u>right</u> end, and in the third case, when the function has both negative and positive values, $S_1$ is <u>inside</u> of the scale (not neccesarily at the halving point). Have a look at the functions e.g. $x^2$, $x^3$, $\sqrt{x}$, $sin(x)$, $log(x)$, $log(sin(x))$, $\sqrt[3]{x^3-1}$, $arctan(x)$ and $tanhyp(x)$. Observe, that the increasing or decreasing property of the function does not affect the questions concerned $x_0$ and $S_1$.

If the function $f$ has no root, the easiest solution we can do is putting $S_1$ <u>anywhere</u> on the scale and measure all (signed) distances from this (geometrical) point. However we are still interested in the relationship of $S_1$ to the function $f(x)$ and we want to find an appropriate <u>symbol</u> (surely not a real number) $x_0$ to write at this point. If some numbers $x_1, x_2, ..., x_n, ...$ are written closer and closer to $S_1$ on the scale, we clearly must have $lim_{n\to\infty} f(x_n) = 0$ ! Now follow the numbers $x_n$ where they "go away": compute $x_0 := lim_{n\to\infty} x_n$ ! Yes, $x_0$ must be this limit, since the numbers $x_n$ approach $x_0$, and their functional values $f(x_n)$, i.e. the distances from $S_1$ are going to vanish.

Summarizing: the symbol $x_0$ must be such that

$$lim_{x\to x_0} f(x) = 0 \qquad (4)$$

would hold. Remark, that monotonic (either in- or decreasing) functions can have at most one (extended) number $x_0$ for which (iv) yields.

For example, for the function $f(x) = 1/x$ we have $x_0="\infty"$ since $lim_{x\to\infty} 1/x = 0$. For $e^x$ and $\Phi(x)$ we have $x_0="-\infty"$ since $lim_{x\to -\infty} e^x = 0$ and $lim_{x\to -\infty} \Phi(x) = 0$. For $x\cdot log(x)$ we have $x_0=0$, which is <u>not</u> an element of Dom, since $lim_{x\to 0} x\cdot log(x) = 0$. Keep in mind, that $S_1$ always has to be present on the scale.

Finally, some functions do not have a solution even of (4), as $f(x) = x+1/x$ or $coshyp(x)$. In these cases we have no choice: $S_1$ must stand alone on the scale, no values $x$ are printed next to it and no $x_0$ is set to it. C'est la vie!

The function $f(x)_{942} = \frac{x^3}{x^2-x-2} = x+1+\frac{3x+2}{x^2-x-2}$, mentioned in (xiii), is an ordinary one: it is strictly monotone increasing on the intervals $x<1-\sqrt{7}$ and $1+\sqrt{7}<x$. Since neither $f_{942}(x_0)=0$ nor $lim_{x\to x_0} f_{942}(x)=0$ in these intervals has a solution for $x_0$, $S_1$ is at the right end of the scale for $x<1-\sqrt{7}$ and at the left end for $1+\sqrt{7}<x$, and the closest values "$1-\sqrt{7}$" and "$1+\sqrt{7}$" have (geometrical) distances from $S_1$ exactly $f_{942}(1-\sqrt{7}) \approx -1.893$ and $f_{942}(1+\sqrt{7}) \approx 6.338$, respectively.

Zooming has also problems with $S_1$ and $x_0$ as we mentioned in Section 2 (o). Since $S_1$ must always be on the scale, we can zoom only a <u>neighbourhood</u> ("close to") of $S_1$ even if we are interested in a part of the scale far from $S_1$. In many cases $S_1$ has the mark $x_0 = \infty$ which is not in our computational range! Even when using a computer ([1] or [6]) and translating the scale, do not forget to zoom all the other scales by the same factor when the slide is in operation (adding or substracting). Moreover, in operation we need to see in fact also $S_1$ - recall the basic operation of slide rules. Try out $10^x$, $arctan(x)+\pi$ or $\Phi(x)$. On the other hand, many functions have both both gently and steeply changing parts which implies rare or dense marks (see (vi) in Section 2. For example, investigate the function $arctan(x)$ on the intervals [-10,+10], [-1000,+1000] and [-10000,+10000]

We mentioned the problem of <u>translating</u> $f(x)$ to $g(x)=f(x)+v$ ($v \in R$) also in Section 2 (o). What is the new $S_1'$ and $x_0'$ for $g(x)$ ? We have to set the (new) distances of the marks "$x$" from $S_1'$ to $d_x'=g(x)=f(x)+v$, while the original distance of "$x$" from $S_1$ was $d_x=f(x)$. This means that we have to choose $S_1'$ in the <u>other</u> direction (of the original scale) by $v$ units. If the (original) scale had some mark "$x_{-v}$" in that other direction by $v$ units, then the new $S_1'$ would fit under this mark "$x_{-v}$". In other words: in <u>this</u> case $x_0' = x_{-v}$.

Recall, that we must have $g(x_{-v})=0$, i.e. equivalently $f(x_{-v}) = -v$. Have a look at the graph of $f(x)$ and observe that the equality $f(x_{-v}) = -v$ does <u>not</u> have any connection to $x_0$ or $v$, it is completely a new equation. In many cases $x_0' = x_{-v} = \pm\infty$ may happen, which causes many problems (as zooming), and does not contradict to the fact that $S_1'$ is at finite (geometrical) distance from $S_1$ and $x_0$, namely at $v$ (units).

The "hardest" case (to understand) is when neither $x_{-v} \in Dom(f)$ nor $x_{-v} = \pm\infty$ does exist satisfying $f(x_{-v}) = -v$. In this case, as we stated in a paragraph above, $S_1'$ will stand alone on the scale of $g(x)$. See again the scales for $1/x$, $1+1/x$, $arctan(x)$ and $arctan(x)+\pi/2$.

All the papers [1] through [6] have problems with $S_1$ and $x_0$, See also the problems and solutions in those papers.

### 4 The scales $-f(x)$ and $f(-x)$

The scale for the function $-f(x)$ would be useful for substraction, e.g. in Newton-Leibniz rule, or for division with the logarithmic scales. Well, by the "inverse use" of the original scale $f(x)$ is also a solution for substracting (see the manual of traditional slide rules). However both for easier use and for better understanding we reveal here the scale for $-f(x)$. The graph of the function $-f(x)$ is the mirror image of the graph of $f(x)$ to the $x$ axis: the number "$x$" is the same but the values $y$ and $-y$ differ in their sign only, and these values $y$ and $-y$ give us the (geometrical) distance of the printed number "$x$" from $S_1$ on the scale.

So, the <u>solution</u> for the scale of $-f(x)$ is the following: just reflect the whole scale of $f(x)$ to $S_1$ and do <u>not</u> alter any printed value "$x$" on the scale. We are lucky to be able to translate scales, otherwise this new scale would hang out of the stator. If $S_1$ was at one end of the scale, after reflecting and translating it will be at the other end. That is, finally we have reflected the scale to the middle point of the old scale to get the new one. Clearly $-f(x)$ is de-/increasing just in $f(x)$ was in-/de- creasing. Try out any function, e.g. $x^3$, $x^\alpha$ for various $\alpha \in R$, $\sqrt[3]{x^3-1}$ and $\sqrt[3]{1-x^3}$, $sin(x)$, $cos(x)$, $arccos(x)$, $arctan(x)$, $\Phi(x)$, $f_{953}(x)=x+\sqrt{1-x^2}$, etc., and their tanslates.

The scales and the graphs of the functions $f(-x)$ and $-f(x)$ must not be confused! However, the construction and the final shape of these scales look very similar.

The graph of the function $f(-x)$ is the mirror image of the graph of $f(x)$ to the $y$ axis: only the <u>sign</u> of the number "$x$" is changed but the value $y$ remain the same. This unchanged value $y$ give us the distance of the printed number "$x$" and "$-x$" from $S_1$ on the scale. So, we get the scale of $f(-x)$ as reflecting the whole scale of $f(x)$ to $S_1$ and <u>altering</u> the signs of all the printed values "$x$" to "$-x$".

From the above paragraphs we can derive that the scales of $-f(x)$ and $f(-x)$ differ <u>only</u> the signs of the printed values "$x$" and "$-x$". The symmetry of the scale can be observed only when there have been already both positive and negative marks on the scale.

The <u>odd</u> functions have the property $-f(x)=f(-x)$, so $f(0)=f(-0)=-f(0)$ which implies $x_0=0$ and so these scales are symmetric to the point $S_1$ (no translation is necessary) since $d_{-x}= f(-x) = -f(x) = -d_x$. The other direction is also true: the scale is symmetric to $S_1$ if and only if $f(x)$ is odd.

Scales for <u>even</u> functions can be constructed for their half domain only since by $f(-x) = f(x)$ their graphs are symmetric to the *y* axis: they can be monotonic of <u>different</u> type for $x<0$ and $0<x$. However, using the equality $f(-x)=f(x)$ we can write (or just imagine) both the $+x$ and $-x$ at the same place (tick) on the scale.

The generalization of odd functions are the functions <u>central symmetric</u> to a point Q($q,r$), that is

$$f(q+z) - r = r - f(q-z), \quad \text{for } z \in \mathbb{R} \tag{5a}$$

or, equivalently

$$f(q-z) - r = r - f(q+z), \quad \text{for } z \in \mathbb{R} \tag{5b}$$

or

$$\frac{f(q-z)+f(q+z)}{2} = r \quad \text{for } z \in \mathbb{R} \tag{5c}$$

and clearly $f(q)=r$. For example $\Phi(x)$ is symmetric to the point Q($0, 1/2$) since $\Phi(-x)=1-\Phi(x)$, $cos(x)$ is symmetric to Q($\pi/2,0$) and so $arccos(x)$ is to Q($0,\pi/2$), $arctan(x)+\pi/2$ is to Q($0,\pi/2$). Odd means to be symmetric to Q($0,0$). Any function $f(x)$ is symmetric to a point Q($q,r$) <u>if and only if</u> the function $h(x):=f(q+x)-r$ is odd, since $h(-x)=-h(x)$ is equivalent to $f(q-x)-r = -(f(q+x)-r) = r - f(q+x)$, as in (5).

Finally, if $f(x)$ is symmetric to Q($q,r$), then its scale is symmetric to the point marked by the value "$q$", since the ticks corresponding to $x_1=q+z$ and to $x_2=q-z$ have distances from the mark "$q$" as

$$d_{x1} = f(q)-f(x_1) = r-f(q+x) = f(q-x)-r = f(x_2)-f(q) = d_{x1}$$

by (5). See e.g. the function $f(x) = 1/100 \cdot (x-2)^3+1$ on Figure 4.

## 5 Special properties

Here we continue our investigations begun in Section 2.

**(xiv)** The (arbitrary base) $log_a(x)$ function's well known property

$$log_a(a^h \cdot x) = h + log_a(x) \quad \text{for all } x,h \in \mathbb{R}, 0<x \tag{6}$$

for integer $h$ makes the "pure" logaritmic scales (named B, K and C) funny and simple. For example the scales for the intervals [1,10] and [10,100] (now $a=10$ and $h=1$) are the "same" (<u>using a translation</u>): only the number of 0-s at the marked numbers are different. This also can be seen on Figure 1: the scale on the lower part of the *y* axis ($0.1 \le x \le 1$) is the same as the upper part ($1 \le x \le 10$). This applies for all $h$, best seen on scale K. Unfortunately no other functions do exist than $f(x)=log_a(x)+c$ satisfying (6). (Each $a$ is convenient to the number system of base $a$.)

**(xv)** The function $f(x)$ is called **homogeneous** of order $k$ ($k \in \mathbb{R}$ fixed) if

$$f(\lambda \cdot x) = \lambda^k \cdot f(x) \quad \text{for all } x,\lambda \in \mathbb{R}, 0<\lambda \tag{7a}$$

or equivalently

$$f(\mu^k \cdot x) = \mu \cdot f(x) \quad \text{for all } x,\mu \in \mathbb{R}, 0<\mu \tag{7b}$$

holds. Clearly then $f(0)=0$, i.e. $x_0=0$. Using (7) again simplifies the construction of the scale for $f(x)$: <u>zooming</u> the scale for $[0,b]$ by the factor $\mu$ ($\mu=10^h$ is adviced) we immediately get the scale for $[0, \mu^k \cdot b]$. If, in the same time we alter the <u>unit</u> $u$ by the reciprocial value $1/\mu$ the original scale returns, only the marks $x$ are changed to $\mu^k \cdot x$. (Like in (xiv) but zooming instead of translation.) The functions $x$, $x^2$, $x^3$, $1/x$ and all the functions $x^\alpha$ for $\alpha \in \mathbb{R}$ and $0 \le x$ satisfy (7) and on their scales the above zooming phenomenon can easily be observed. Unfortunately no other functions do exist than $f(x)=c \cdot x^\alpha$ satisfying (7).

This zooming property was a withdrawal in our constructions in [3], that is why we had to look for "more irrational" numbers instead of "easy rational" ones.

**(xvi)** Some other functions have the property

$$f(\lambda+x) = h_\lambda \cdot f(x) \quad \text{for all } x,\lambda \in \mathbb{R} \tag{8}$$

where the constant $h_\lambda$ depends only on $\lambda$ but not on $x$. This equality implies that all intervals $[v+\lambda,w+\lambda]$ are zoomed images of $[v,w]$ for all $v,w,\lambda$ by the factor $h_\lambda$. For example, the functions $a^x$ satisfy (8) with $h_\lambda=a^\lambda$. Since $x_0 = $ "$-\infty$" for these functions, the <u>equal length</u> scales for $[x_0, x_1]$ and $[x_0, x_1+\lambda]$ are exactly the same, since for achieving the equal lengths of the scales we had to decrease the unit $u$ by the factor $1/h_\lambda$ (as in (xv)). This property can be examined by [6], see also Figure 2. Unfortunately no other functions do exist than $f(x)=c \cdot a^x$ satisfying (8).

**(xvii)** We discussed the role of the derivative $f\,'$ in entries (vi) through (viii). $f(x)$ is **convex/concave** ("watching from below") if its first derivative $f\,'$ is increasing/decreasing, so the <u>denseness</u> of the numbers (or ticks) on the scale is increasing/decreasing. An illustrative example is $f_{921}(x)=x^5-3x^3$ for $-3/\sqrt{5}<x<3/\sqrt{5} \approx 2.236$, which is monotone decreasing and has both convex and concave parts on this interval.

## 6 The inverse function $f^{-1}$

Inverting means changing $x$ and $y$ as $y=f(x) \Leftrightarrow f^{-1}(y)=x$. (Well, after solving this equality we usually write $f^{-1}(x)$ instead of $f^{-1}(y)$.) Each strictly monotone function has a unique inverse and of the same type (in- or decreasing). Practical investigations show that the density of marks (ticks) changes: the scale of $f(x)$ is dense/rare if and only if the scale of $f^{-1}(x)$ is rare/dense. It is right by the theory, too, since the derivative of the inverse is just the reciprocal of the original function's: $(f^{-1})' = 1/f'$ and see (vi) in Section 2.

Is there a method to obtain the scale of $f^{-1}(x)$ from the "original" scale (of $f(x)$) ? The construction of the <u>graph</u> of $f^{-1}(x)$ is simple: reflect to the straight line "$y=x$", i.e. looking the paper from the <u>back</u> side and rotating the axes $y$ and $x$. The <u>scale</u> for $f^{-1}(x)$ could be constructed on the base of Figure 1 as: start drawing horitontal lines from the equidistant numbers 0,1,2,... on the $y$ axis to the graph of the function $f(x)$, then the vertical lines from these points meet the scale for $f^{-1}(x)$ on the $x$ axis.

The "<u>RubberBand</u>" method explains the properties of the <u>scale</u> for $f^{-1}(x)$ from other point of view. We see only the <u>numbers</u> "$x$" on the scale for $f(x)$ but not the real, geometrical distance scale 0,1,2,... (in some units). By a pencil sketch this equidistant (geometrical) scale under the scale of $f(x)$. We have to reverse the equality $d=d_x=y=f(x)$ to $x=f^{-1}(y) = f^{-1}(d)$. Think on the scales as they were printed on (real physical) rubber band and pull-push it in such a way until the marks (ticks) "$x$" on the scale for $f(x)$ form an equidistant scale, as in everyday measure tape. Now the other scale, containing the pencil marks "$d$" will show you the scale for the function $f^{-1}(x)$.

The hard question is to find $S_1'$ and $x_0'$ for the inverse scale. By definition $S_1'$ is the (geometrical) point where $x_0'$ is marked on the new scale with the property $f^{-1}(x_0')=0$, that is $x_0'=f(0)$ assuming $0 \in \text{Dom}(f)$. In this case $S_1'$ belongs to the tick denoted by 0 on the original scale (for $f(x)$) and the value is $x_0'=f(0)$. (Keep in mind that all ticks are moved somewhere, according to the previous paragraph.) In the case $0 \notin \text{Dom}(f)$ we can only use the equality

$$\lim_{x \to a} f(x) = b \quad \Leftrightarrow \quad \lim_{x \to b} f^{-1}(x) = a \qquad \text{for } a,b \in \mathbb{R} \cup \{\pm\infty\} \tag{9}$$

for determining $S_1'$ and $x_0'$. See e.g. the functions $lg(x)$, $1/x$ or $\sqrt{x-1}$.

The *self inverting* functions (i.e. $\forall x\, f(x)=f^{-1}(x)$) have interesting properties: they are symmetric to the straight line "$y=x$". The only increasing self inverting function is $f(x)=x$ (why?). Some self inverting functions are: $t/x$ and $\sqrt[\alpha]{t-x^\alpha}$ for any fixed $\alpha, t \in \mathbb{R}$. Clearly the scales for $f^{-1}(x)$ and $f(x)$ must be the same for self inverting functions, but it is worth to have a look at the scale of $\sqrt[3]{1-x^3}$ on Figure 3. The dense (at $x=0$) and the rare (at $x=1=S_1$) ticks fields are interchanged, but because of the "RubberBand" method the denseness around these ticks are interchanged and set everything back to the original. Each self inverting function satisfies (no $f^{-1}$ in (10)):

$$f(a)=b \quad \Leftrightarrow \quad f(b)=a \qquad \text{for } a,b \in \mathbb{R} \tag{10}$$

## 7 Applications

For what to use these (infinitely many) new scales? Wait for engineers to order new scales for their own?

We have a collection in [2] for several calculation usage: the $x^\alpha$ scales for many mathematical, statistical, physical and other applications (springs and resistors $1/D_3=1/D_1+1/D_2$ , inductivities $\sqrt{L_3}=\sqrt{L_1}+\sqrt{L_2}$ , Pithagoras rule, $D^2(\xi_1+\xi_2) = D^2(\xi_1)+D^2(\xi_2)$ in statistics). The (real) $lg(sin(x))$ scale is good for the sine law (of triangles) or Snellius-Descartes's law. The cumulative distribution functions $\Phi(x)$ and $F_\lambda(x)=1-e^{-\lambda x}$ and others in probability theory for the Newton-Leibniz rule. See also [4] and [6].

However many new scales are useful for <u>illustrating</u> formulas we learn in secondary or high school level, beyond the above mentioned ones. Though this question is close to the theme "analog vs. digital" we do not want to discuss in detail here. Only some examples, without order or completeness.

 **a)** use the scales $sin^2(x)$ and $cos^2(x)$ to illustrate the equality $sin^2(x)+cos^2(x)=1$, or the scales $sinhyp^2(x)$ and $coshyp^2(x)$ to $sinhyp^2(x)+1=coshyp^2(x)$.

 **b)** use the scales $sin^2(x)$ and $cos(2x)/2$ to illustrate the equality $sin^2(x)=(1-cos(2x))/2$, similar to the hyperbolic functions,

 **c)** use the scale $\Phi(x)$ to observe $\Phi(-x)=1-\Phi(x)$, and for illustrating any odd or symmetric to a point functions (see at the end of Section 4),

 **d)** in general: every question (i) through (xvii) illustrates in a new form the basic properties of functions,

**e)** the scales $lg(x)$ and $1/x$ and every $x^\alpha$ have fractal-like properties: zoomings of them result similar scales (see Section 5).

## 8  Some examples

For illustrating the arguments in this article we show here some example scales, generated by [6]. Let us emphasize that [6] is a pure program (application) for mathematical research, while Hoffman's [1] program is a professional one for demonstrating traditional and some new scales.

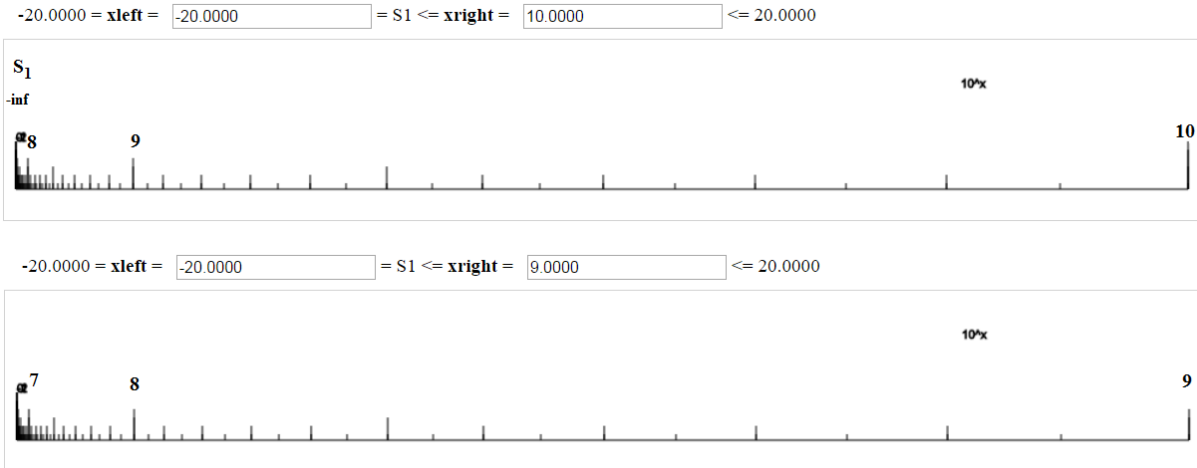

**Figure 2.:** The $10^x$ scales for the intervals [-∞,10] and [-∞,9]

No $x_0$ such that $f(x_0)=0$ would hold for $f(x)=10^x$ but $lim_{x\to-\infty}=0$ so "-inf" is written at $S_1$. The density of the numbers exploids when going to "-inf" since the functions is almost horizontal  (has a horizontal asymptote). The below scale is a zoomed version of the upper one but looks similar, illustrating (8) in Section 5, using "-∞+λ=-∞ for any λ∈R" .

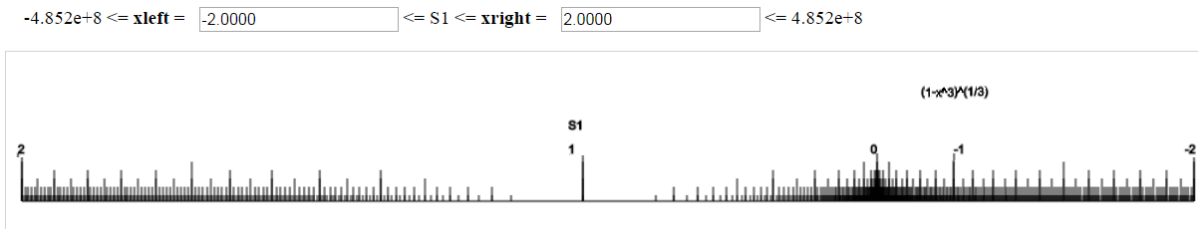

**Figure 3.:** $\sqrt[3]{1-x^3}$  for  -2≤x≤+2

$x_0=1$ since $f(1)=0$ for $f(x)=\sqrt[3]{1-x^3}$ and $S_1$ is in the middle of the scale since $f(x)$ had both negative and positive values.  $f(x)$ is <u>decreasing</u> so the values ("the remainder of the $x$ axis") are increasing from right to left! The values and ticks around $x=1$ are rare since $f(x)$ is very steep here (in fact $f\,'(1)=+\infty$) while around $x=0$ the ticks are crowded since $f(x)$ is almost horizontal here (in fact $f\,'(0)=0$).  Far from 1 and -1 the scale is almost equidistant since $f(x)$ has the slanted asymptote  $f(x) \approx -x$ .

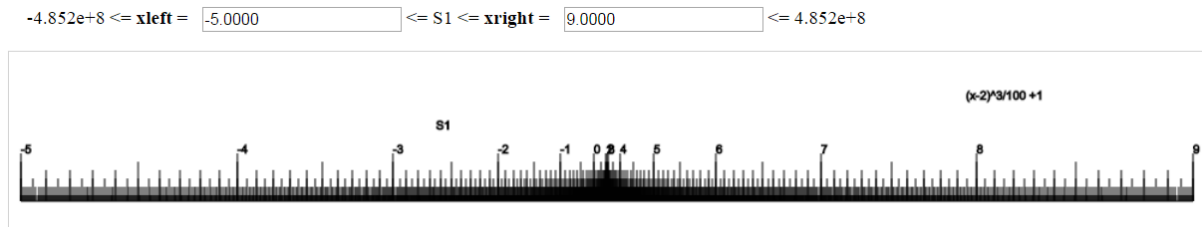

**Figure 4.:** $1/100 \cdot (x-2)^3 + 1$  for  $-5 \leq x \leq +9$

The function $f(x) = {}^1\!/_{100} \cdot (x-2)^3 + 1$ has a horizontal tangent line at $x=2$ that is why the ticks are crowded at this point. When going far away the numbers are getting rare since the slope is going to $\pm\infty$. $S_1$ is at the point $x_0 = 2-\sqrt[3]{100} \approx -2.6416$ since $f(x_0)=0$. Finally, $f(x)$ is symmetric to the point Q(2,1) as the scale is symmetric to $x=2$ and $f(2)=1$.

# 9 Acknowledgement

**Notes**

(*) University of Pannonia, Veszprém, Hungary, szalkai@almos.uni-pannon.hu

(**) For any function $f$ we denote the **inverse** of $f$ by $f^{-1}$ which means $f^{-1}(y)=x \leftrightarrow f^{-1}(y)=x$ for any $x \in \text{Dom}(f)$ and $y \in \text{Im}(f)$.

(***) $f(x) \approx cx+b$ means $\lim_{x \to \infty} f(x)/x = c$ and $\lim_{x \to \infty} (f(x)-c \cdot x) = b$.